\documentclass[a4paper]{article}
\makeatletter{}\usepackage			{amsmath}
\usepackage			{amsfonts}
\usepackage			{array}
\usepackage			{comment}
\usepackage			{enumerate}
\usepackage[mathscr]		{eucal}
\usepackage			{graphicx}
\usepackage[utf8]		{inputenc}
\usepackage			{latexsym}
\usepackage			{longtable}
\usepackage			{multirow}
\usepackage			{rcs}
\usepackage			{theorem}
\usepackage			{xspace}
\usepackage[lined,boxed]	{algorithm2e}
\usepackage			{url}

\newtheorem{theorem}{Theorem}

\newtheorem{definition}{Definition}

\DeclareMathOperator*	{\epi}	{{\mathrm epi}}

\newcommand	{\xs}		{{x^\star}}

\newcommand	{\epifs}	{{\epi f^\star}}

\usepackage{xspace}

\makeatletter{}

\newtheorem{lemma}{Lemma}

\newcommand{\EndProof}{\rule{2ex}{2ex}\xspace}
\newcommand{\bleq}[2]{\begin{equation}\label{#1}{#2}\end{equation}}

\DeclareMathOperator	{\Dom}	{{\mathrm dom}}
\DeclareMathOperator	{\intr}	{{\mathrm int}}

\newcommand{\R}{\mathbb{R}}

\newcommand{\Dfrac}[2]{{\displaystyle\frac{{#1}}{{#2}}}}

\newcommand{\fs}{{f^\star}}

\newcommand{\dep}{{\partial}}

\newcommand{\Proof}{{\bf Proof.}\xspace}

\author{E.A. Nurminski
\thanks{Far Eastern Federal University, School of Natural Sciences,
Ajax St., Vladivostok, Russky Island, Russia}}
\title{A Conceptual Conjugate Epi-Projection Algorithm of Convex Optimization:
Superlinear, Quadratic and Finite Convergence
\thanks{This work is supported by RF Ministry of Education and Science,
project 1.7658.2017/6.7}
}
\begin{document}
\maketitle
\begin{abstract}
\makeatletter{}This paper considers a conceptual version of a convex optimization algorithm
which is based on replacing a convex optimization problem
with the root-finding problem for the approximate sub-differential mapping
which is solved by repeated projection onto the epigraph of conjugate function.
Whilst the projection problem is not exactly solvable in finite space-time
it can be approximately solved up to arbitrary precision by simple iterative
methods, which use linear support functions of the epigraph.
It seems therefore useful to study computational characteristics of the idealized
version of this algorithm when projection on the epigraph is computed precisely to estimate the potential benefits
for such development.
The key results of this study are that the conceptual algorithm
attains super-linear rate of convergence in general convex case,
the rate of convergence becomes quadratic for objective functions
forming super-set of strongly convex functions,
and convergence is finite when objective function has sharp minimum.
In all cases convergence is global and does not require differentiability of the objective.
 
 
\end{abstract}
\makeatletter{}\section*{Introduction}
\label{intro}
We consider a finite-dimensional nondifferentiable convex optimization problem (COP)
\bleq{copt}{ \min_{x\in E} f(x) = f_\star = f(\xs), \xs \in X_\star\,,}
where
\(E\) denotes a finite-dimensional space of primal variables and \(f:E\to\R\) is a finite convex function,
not necessarily differentiable.
As we are interested in computational issues related to solving (\ref{copt})
mainly we assume that this problem is solvable and has nonempty and bounded set of solutions \(X_\star\).

This problem enjoys a considerable popularity due to its important theoretical properties and numerous
applications in large-scale structured optimization, discrete optimization,
exact penalization in constrained optimization, and others.
This led to the development of different algorithmic ideas, starting with the subgradient algorithm due
to Shor (see \cite{Shor2012} for the overview and references to earliest works)
and followed by
conjugate subgradient algorithms \cite{wolfe,Li2013}
bundle methods \cite{LemUr},
space dilatation and \(r\)-algorithms
\cite{ndo2003},
\(\epsilon\)-subgradient methods \cite{rzev, kiwiel, nurmi91},
\(VU\)-methods
\cite{Miff+Saga2005},
proximal point algorithms
\cite{Parikh+Boyd-2013}
and many others.
These algorithms were widely used for solving many academical and practical problems,
however only in a few cases the estimates for the rate of convergence were obtained.
The most notable case is probably algorithms of proximal point family (PPA),
which use the smooth approximation of the original COP by Moreau-Yosida regularization.
In this case the superlinear rate of convergence was attained both
for conceptional and implementable versions of PPA
\cite{burke+Qian2000}.

In this paper we suggest
another algorithm with attractive rate of convergence at the conceptual level.
It was suggested by the author in \cite{nurmi91,nurmi97},
but at that time only superlinear convergence for the general convex case was established.
The finite convergence was demonstrated for the conical polyhedral functions only.
The main attention was given to the implementable versions of the algorithm
which used polytope approximations of the epigraph of the conjugate function.
The following developments
\cite{vor14, vor16}
of the implementable versions
concentrated on different computational details
and applications of these algorithms to some specific problems.
However, the positive computational experience obtained,
and some new computational ideas
\cite{vornur15,nurmi17}
which introduced additional cuts into approximation schemes for
the epigraph of the conjugate functions, revived the interest into
the conceptual version of the main algorithmic idea.
The new approach to the study of the conjugate epi-projection algorithm
resulted in the new generic estimate of the algorithm progress toward the solution,
formulated in terms of directional derivatives of the conjugate function.
The new analysis, based on this estimate, reported here,
showed that the algorithm  does not only attain superlinear rate of convergence in
quite general convex case,
but
the rate of convergence becomes quadratic for objective functions which are slightly
more general than strictly convex in a vicinity of optimal solution,
and the convergence is finite when objective has sharp minimum.
In all cases the convergence is global and does not require differentiability of the objective.
\makeatletter{}\section{Notations and Preliminaries}
\label{nots}

Throughout the paper we use the following notations:
\(E\) is a finite dimensional euclidean space of primal variables
of any dimensionality.
The inner product of vectors \(x, y\) from \(E\) is denoted as \(xy\).
The cone of non-negative vectors of \(E\) is denoted as \(E_+\).
The set of real numbers is denoted as \(\R\) and \(\R_\infty = \R \cup \{\infty\} \).

The norm in \(E\) is defined in a standard way: \(\|x\| = \sqrt{xx}\)
and for \(X \subset E\) \(\|X\| = \sup_{x\in X} \|x\|\).
This norm defines of course the standard topology on \(E\) with the common
definitions of open and closed sets and closure and interior of subsets of \(E\). 
The interior of a set \(X\) is denoted as \(\intr(X)\).

The unit ball in \(E\) is denoted as \(B = \{x: \|x\| \leq 1 \}\).
The support function of a set \(Z \subset E \) is denoted and defined as
\( (Z)_x = \sup_{z \in Z} xz \).

A vector of ones of a suitable dimensionality is denoted by \( e = ( 1,1, \dots, 1).\)
A standard simplex \(\{ x: x \geq 0, xe = 1 \}\) with \(x \in E, \dim(E) = n\)
is denoted by \( \Delta_E. \)

We use the standard definitions of convex analysis
(see f.i. \cite{LemUr})
related mainly to functions \(f: E \to \R_\infty\):
the domain \(\Dom f \) of a function \(f\) is the set \(\Dom f = \{ x: f(x) < \infty \}\),
the epigraph \(\epi f\) of a function \(f\) is a set
\( \epi f = \{(\mu, x): \mu \geq f(x) \} \subset \R_\infty \times E\).

Further on all functions are convex in a sense that their epigraphs are {\em convex} subsets
of \(\R_\infty \times E\).
\begin{definition}
For a convex function \(f: E \to \R\) and fixed \(x \in E\) the set
\( \dep f(x) = \{ g: f(y) - f(x) \geq g(y-x) \mbox{ for all } y \in \Dom f \}\)
is called a sub-differential of \(f\) at the point \(x\). 
\label{defsubg}
\end{definition}

The sub-differential of \(f\) is well-defined and is a closed bounded convex set for all \(x \in \intr(\Dom f) \).
At the boundary of \(\Dom f\) it may or may not exists.
The sub-differential of \(f\) is also upper semi-continuous as a multi-function of \(x\) when exists.

\begin{definition}
The directional derivative of
a finite convex function \(f\) 
at point \(x\) in direction \(d\)
is denoted and defined as 
\[
\dep f(x; d) = \lim_{\delta \to +0} (f(x + \delta d) - f(x))/\delta.
\]
\label{dirder}
\end{definition}
It is known from convex analysis that
\( \dep f(x; d) = \sup_{g \in \dep f(x)} gd  = (\dep f(x))_d \).

\begin{definition}
For a convex function \(f: E \to \R_\infty\) the function
\bleq{cf}{ \fs(g) = \sup_x \{ gx - f(x) \} = (\epi f)_{\bar g}, \mbox{ where } \bar g = (-1, g) \in \R_\infty \times E }
is called a conjugate function of \(f\). 
\label{defconf}
\end{definition}

The key result of convex analysis is that for a closed function \(f\)
which epigraph \(\epi f\) is a closed set
\bleq{duality}{ \sup_g \{ gx - \fs(g) \} = (\epi \fs)_{\bar x} = f(x), } 
where \(\bar x = (-1, x) \in \R_\infty \times E \).

It is also easy to see that if \( (\epi \fs)_{\bar x} = g_x x - \fs(g_x) \) then \(g_x \in \dep f(x)\)
and the other way around:
for \(\bar g = (-1, g)\)
if \( (\epi f)_{\bar g} = g x_g - f(x_g) \) then \(x_g \in \dep \fs(g)\).

The trivial consequence of the Definition \ref{defconf} is that \( \fs(0) = -\inf_x f(x) \)
 which is the key correspondence used by the conjugate epi-projection algorithm, considered further on.
As the conjugate epi-projection algorithm operates in the conjugate space its convergence properties depend upon
the properties of the conjugate function of the objective.
Therefore we introduce some additional classes of primal functions to ensure the desired behavior
of the conjugates.
 
\begin{definition}
\label{sup-def}
Convex function \(f\) is called sup-quadratic with respect to a point \(x \in \intr(\Dom f)\)
if there exists a constant \(\tau > 0\) such that 
\bleq{sup-q}{f(y) - f(x) \geq g(y-x) + \frac12 \tau \|y - x\|^2}
for any \(g \in \partial f(x)\) and any \(y\).
\end{definition}
We will call \(\tau\) the sup-quadratic characteristic of \(f\) at \(x\).
Notice that strongly convex functions are sup-quadratic at any \(x\)
from their domains, however a function \(f\), sup-quadratic at some \(x\), need not to be strongly convex.

A symmetric definition can be given for {\em sub-quadratic} functions.
\begin{definition}
\label{sub-def}
Convex function \(f\) is called sub-quadratic with respect to a point \(x \in \intr(\Dom f)\)
if there exists a constant \(\tau > 0\) such that 
\bleq{sub-q}{f(y) - f(x) \leq g(y-x) + \frac12 \tau^{-1}\|y - x\|^2}
for any \(y \in \Dom f \) and some \(g \in \partial f (x)\).
\end{definition}
Notice that it follows from this definition that
the function \(f\),
sub-quadratic at point \(x\)
is in fact differentiable at this point.
Of course not all functions differentiable at \(x\) are
sub-quadratic.

Definitions \ref{sup-def} and \ref{sub-def}
allow us to establish an important properties of conjugates
functions for sup-quadratic primals.
\begin{lemma}
\label{sup-sub}
Let \(f: E \to \R \) attains its minimum value \(f_\star\) at the point \(\xs\) 
and \(f\) is sup-quadratic at point \(\xs\) with the positive sup-quadratic characteristic \(\tau\).
Then \(\fs(g)\) is sub-quadratic at \(g = 0\) with the corresponding
sub-quadratic characteristic not lower then \(\tau^{-1}\).
\end{lemma}
\Proof
By definition for any \(x\)
\bleq{gt}{
\frac12 \tau \|\xs - x\|^2 \leq
f(x) - f_\star = f(x) + \fs(0)
} 
and hence
\bleq{chin}{
\fs(g) - \fs(0) = x_g g - (f(x_g) + \fs(0)) \leq x_g g - \frac12 \tau \|\xs - x_g\|^2
}
for any \(x_g \in \partial f^\star(g)\).
Hence
\bleq{chin2}{ 
\begin{array}{c}
\fs(g) - \fs(0) \leq
\xs g + (x_g - \xs) g - \frac12 \tau \|\xs - x_g\|^2 \leq
\\
\xs g + \sup_z \{ z g - \frac12 \tau \|z\|^2 \} =
\xs g + \frac12 \tau^{-1} \|g\|^2.
\end{array}
}
\EndProof

Another interesting subclass of convex functions are those which have zero in the interior
of the subdifferential at the solution \(\xs\) of a COP (\ref{copt}), that is
\( 0 \in \intr( \partial f(\xs)). \)
This condition is also known as ''sharp minimum''
and extended further on in
\cite{Ferris88} and others.
The special attraction of this case is that the known proximal method has then a finite termination 
\cite{Ferris1991} for such problems.

We notice now that the conjugate functions for objectives with sharp minimum have very simple behavior
in the vicinity of zero which also guarantees the finite termination of
the conjugate epi-projection algorithm as well.
\begin{lemma}
\label{lincon}
If solution \(\xs\) of (\ref{copt}) is such that \( 0 \in \intr( \partial f(\xs)) \)
then there is \(\rho > 0\) such that \(\fs(g) = \sup_x \{ gx - f(x) \} = g \xs - f(\xs)\)
for \(\|g\| < \rho\).
\end{lemma}
\Proof
If \(\rho\) is small enough then sharp minimum condition implies
\(0 \in \partial (f(\xs) - g\xs) = \partial f(\xs) - g \) for any \(g \in \rho B\) and therefore
\[
\fs(g) = \sup_x \{ gx - f(x) \} = g \xs - f(\xs)
\]
is a linear function of \(g\).  
\EndProof

For additional results on connections between sharp minimum and properties of conjugate functions see also \cite{Zhou2012}.
 
\makeatletter{}\section{Conjugate Epi-Projection Algorithm}
\label{conspa}

\newcommand{\xikp}{\xi_k^p}
\newcommand{\gpk}{g_p^k}
\newcommand{\xpk}{x_p^k}
\newcommand{\zk}{z^k}

As it was already mentioned the basic idea of the 
conjugate epi-projection algorithm
consists
in considering the convex problem (\ref{copt})
as the problem of computing the conjugate function of the objective
at the origin:
\[
\fs(0) = -\min_x f(x) = -f_\star = \inf_{(0, \mu) \in \epifs} \mu.
\]
We suggest to use for computing \(\fs(0)\)
the algorithm based on projection onto
the epigraph
\(\epi \fs\).
This idea 
demonstrates some promises for effective solution of (\ref{copt}) and
suggests some new computational ideas.

This version of the algorithm consists in execution of an infinite sequence of iterations,
which generates the corresponding sequence of points \(\{ (\xi_k, 0) \in \R\times E, k=0,1, \dots\}\)
with \(\xi_k \to \fs(0)\) when \(k\to\infty\).
For each of these iterations it calls a subgradient oracle which for any \(x \in E\)
computes \(f(x)\) and arbitrary \(g \in \partial f(x)\).
Also it requires solution of nonlinear projection problem which makes the algorithm
strictly speaking unimplementable.
However the analysis of the algorithm demonstrate its potential and can show the ways to
its practical implementations.
The principal details of the iteration of the conjugate epi-projection algorithm are given on the Fig. Algorithm \ref{cco-algo}.
\begin{algorithm}
\SetAlgoLined
\KwData{The convex function \(f:E \to \R \), the epigraph \(\epifs\),
the current iteration number \(k\) and the current approximation \(\xi_k \leq \fs(0)\).}
\KwResult{The next approximation \( \xi_{k+1} \) such that \( \xi_k \leq \xi_{k+1} \leq \fs(0)\)
}
\par\noindent
Each iteration consists of two basic operations: {\bf Project} and {\bf Support-Update}

{\bf Project.} 
Solve the projection problem of the point \((\xi_k, 0)\) onto \(\epifs\):
\[
\min_{(\xi, g) \in \epifs} \{(\xi - \xi_k)^2 + \|g\|^2 \} = (\xikp - \xi_k)^2 + \|\gpk\|^2
\]
with the corresponding solution \((\xikp, \gpk) = (\fs( \gpk), \gpk) \in \epifs\).
We demonstrate in the analysis of the algorithm convergence that \(\fs(0) \geq \xikp > \xi_k\)
if \(\xi_k < \fs(0)\).

{\bf Support-Update}
Compute support function of \(\epifs\) with the support vector
\( z^k = -(\xikp - \xi_k, \gpk) \in \R \times E\)
\[
\begin{array}{c}
(\epifs)_{z^k} = \sup_{(\mu, g) \in \epi\fs} \{ -(\xikp - \xi_k) \mu + \gpk g) \}
= 
\\
(\xikp - \xi_k) \sup_{(\mu, g) \in \epi\fs} \{ - \mu + \Dfrac{\gpk}{(\xikp - \xi_k)} g \}
= 
(\xikp - \xi_k) \sup_{(\mu, g) \in \epi\fs} \{ - \mu + x_p^k g \}
=
\\
(\xikp - \xi_k) ( x_p^k \tilde \gpk - \fs(\gpk)\} = (\xikp - \xi_k) f(\xpk),
\end{array}
\]
where \(\xpk = \gpk/(\xikp - \xi_k)\).
Notice that as \(f\) is assumed to be a finite function this operation is well-defined.

Finally we update the approximate solution with \(\xi_{k+1}\) using the relationship
\[
\bar\xi_{k+1} z^k = (\epifs)_{z^k}, \mbox{ where } \bar\xi_{k+1} = (\xi_{k+1}, 0) \in \R\times E,
\]
which actually amounts to \(\xi_{k+1} = -f(\xpk)\),
increment iteration counter \(k \to k+1\), etc.
\caption{The basic iteration of the conceptual conjugate epi-projection algorithm algorithm}
\label{cco-algo}
\end{algorithm}
For better understanding of these two operations
they are illustrated on the Fig \ref{cnc:fig1}, \ref{cnc:fig2}.
\begin{figure}
\centering
\parbox{.85\textwidth}{\begin{minipage}[t]{.4\textwidth}
\includegraphics[scale=0.35]{concept-2.mps}
\caption{{\bf Projection.}
Solution of projection problem
\( \min_{(\xi, g) \in \epifs} \{(\xi - \xi_k)^2 + \|g\|^2 \}.  \) 
}
\label{cnc:fig1}
\end{minipage}\qquad
\begin{minipage}[t]{.4\textwidth}
\includegraphics[scale=0.35]{concept-3.mps}
\caption{{\bf Support-Update.}
Compute support function \( v_k = (\epifs)_{z^k}\) and
update the approximate solution with \(\xi_{k+1}\)
\( \xi_{k+1}  = v_k/(\fs(\gpk) - \xi_k).  \)
}
\label{cnc:fig2}
\end{minipage}
}
\end{figure}

Convergence of the Algorithm \ref{cco-algo} is confirmed by the following theorem.
\begin{theorem}
Let \(f\) be a finite convex function with the finite minimum
\( f_\star = \min_x f(x) = -\fs(0)\) and \(\xi_k, k=1,2, \dots\) are defined by the Algorithm \ref{cco-algo}
with \(\xi_0 < \fs(0)\).
Then 
\(\lim_{k\to\infty} \xi_k = \fs(0) = -f_\star\).
\label{thm-algo}
\end{theorem}
\Proof
Assume that on \(k\)-th iteration we have \(\xi_k < \fs(0)\) as the approximation of \(\fs(0)\).
According to Algorithm 1 to construct the next (\(k+1\)-th) approximation \(\xi_{k+1}\)
the point \((\xi_k, 0) \in \R \times E\) is to be projected onto \(\epifs\) first:
\bleq{k-proj}{
\min_{(\xi, g) \in \epifs} \{ (\xi -\xi_k)^2 + \|g\|^2 \} = (\xikp - \xi_k)^2 + \|\gpk\|^2
}
As a result the auxiliary point \((\xikp, \gpk) = (\fs(\gpk, \gpk) \in \epifs \)
is obtained which satisfies optimality conditions
\bleq{o-proj}{
(\fs ( \gpk) - \xi_k ) (\xi - \xikp) + \gpk( g - \gpk) \geq 0
}
for any \((\xi, g) \in \epifs\). 

It is easy to see that \(\xikp > \xi_k\).
Indeed the opposite strict inequality \(\xikp < \xi_k\) contradicts the optimality of \((\xikp, \gpk)\)
as in this case \((\xi_k, \gpk) = (\xikp + (\xi_k - \xikp), \gpk) \in \epifs \),
and
\[
(\xi_k - \xi_k)^2 + \|\gpk\|^2 <
(\xi_k - \xikp)^2 + \|\gpk\|^2 =
\min_{(\xi, g) \in \epifs} \{(\xi_k - \xi)^2 + \|g\|^2 \}.
\]
If \(\xikp = \xi_k\) then \(\R \times \{0\}\) is strictly separable from \(\epifs\):
\[
\xi ( \xi_k - \xikp ) + 0 \gpk = 0 <  \| \gpk \|^2 \leq
\mu ( \xi_k - \xikp ) + g \gpk
\]
for any \((\mu, g) \in \epifs \) as it follows from projection conditions.
Hence \(0 \notin \Dom(\fs)\) which contradicts the assumptions of the theorem. 

According to Algorithm
\ref{cco-algo}
the next approximation \(\xi_{k+1}\) is determined from the equality
\[
(\xikp - \xi_k)(\xi_{k+1} - \xi_k)) - \|\gpk\|^2 = (\xikp) - \xi_k)^2 + \|\gpk\|^2
\]
which gives the following expression for \(\xi_{k+1}\):
\[
\xi_{k+1} = \xi_k + \|\gpk\|^2/(\xikp - \xi_k) \geq \xi_k,
\]
and \(\xi_{k+1} = \xi_k\) if and only if \(\gpk = 0\) which means that we already obtained the solution.

Repeating this operation we obtain the monotone sequence \(\xi_k, k = 0, 1, \dots\) such that
\[
\xi_k \leq \xi_{k+1} \leq \fs(0), k = 0,1, \dots
\]
where inequalities turn into equalities only if either \(\xi_k = \fs(0)\) or \(\xi_{k+1} = \fs(0)\)
which of course makes no difference.
Under these conditions \(\lim_{k\to\infty} \xi_k = \fs(0)\)
which proves the convergence of the algorithm \ref{cco-algo}.
\EndProof

Theorem \ref{thm-algo}
established the convergence of the Algorithm \ref{cco-algo}
under very general conditions,
however to estimate the rates of convergence we need to derive
more convenient estimates for decrease of convergence indicators.
This is provided by the following lemma.
\begin{lemma}
\label{estilem}
Let all assumptions of the theorem \ref{thm-algo} be satisfied 
and \(\xi_k, k=1,2, \dots\) are defined by the Algorithm \ref{cco-algo}
with \(\xi_0 < -f_\star \).
Then 
\bleq{desti}{
\fs(0) - \xi_{k+1} \leq \|\gpk\|(\partial \fs(\gpk; \zk) - \partial \fs(0; \zk)), k = 1,2, \dots
}
where \(\zk = \gpk/\|\gpk\|\).
\label{thm-superlin}
\end{lemma}
\Proof
By construction
\(
\xi_{k+1} = -f(\xpk) = \fs(\gpk) - \xpk \gpk,
\)
where \(\xpk = -\gpk/(\xikp - \xi_k) \).
Then
\[
\fs(0) - \xi_{k+1} = \fs(0) - \fs(\gpk) +  \xpk \gpk \leq \xpk \gpk - \xs \gpk
\]
for any \(\xs \in \partial \fs(0) = X_\star \).
Taking infinum of the right hand side with respect to \(\xs \in \partial \fs(0) \) obtain
\bleq{d-d}{
\begin{array}{c}
\fs(0) - \xi_{k+1} \leq \xpk \gpk - \partial \fs(0; \gpk) \leq  
\sup_{x \in \partial \fs (\gpk)} x \gpk - \partial \fs(0; \gpk) =  
\\
\partial \fs (\gpk; \gpk) - \partial \fs(0, \gpk) =
\|\gpk\|(\partial \fs(\gpk; \zk) - \partial \fs(0; \zk), k = 1,2, \dots,
\end{array}
}
where \(\xpk = -\gpk/(\xikp - \xi_k) \)
and where we used linear positive homogeneity of \( \partial \fs (\cdot; \cdot) \)
with respect to its second argument.
\EndProof

The inequality (\ref{desti}) can be rewritten as
\bleq{conrate}{
\fs(0) - \xi_{k+1} \leq \| \gpk \| (\partial \fs (\gpk; \zk) - \partial \fs(0; \zk))
= \| \gpk \| \theta(\gpk; \zk),
}
and depending on properties of the accuracy multiplicator \(\theta(\gpk; \zk)\)
the convergence rates of Algorithm \ref{cco-algo}
will have different estimates.

First we establish super-linear rate of convergence of Algorithm \ref{cco-algo} for the most general case of a finite
objective function \(f\).
\begin{theorem}
Let all assumptions of the theorem \ref{thm-algo} be satisfied 
and \(\xi_k, k=1,2, \dots\) are defined by the Algorithm \ref{cco-algo}
with \(\xi_0 < -f_\star \).
Then 
\( \fs(0) - \xi_{k+1} \leq \lambda_k (\fs(0) - \xi_k) \) with \(\lambda_k \to 0\) when \(k \to\infty\).
\label{thm-super-linear}
\end{theorem}
\Proof
For the finite \(f\) and bounded nonempty \(X_\star\) in the problem (\ref{copt})
the conjugate function \(\fs\) has nonempty 
\(\Dom(\fs)\) and \(0 \in \intr(\Dom(\fs))\).

Then due to convergence of Algorithm \ref{cco-algo} \( \gpk \to 0 \) when \(k \to \infty\).
In the notations of Algorithm \ref{cco-algo} 
\( \fs(2\gpk) - \fs(\gpk) \geq \fs(\gpk) - \fs(0) \)
by convexity
and hence 
\[
\gpk \xs \leq \fs(\gpk) - \fs(0) \leq \fs(2\gpk) - \fs(\gpk) \leq p^k \gpk
\]
for any \(\xs \in \partial \fs(0) \) and \(p^k \in \partial \fs(2\gpk) \).

After division by \(\|\gpk\| > 0\) it gives
\bleq{left-right}{ \zk \xs  \leq p^k \zk }
where \(z^k = \gpk/\|\gpk\|\).

Taking supremum of the left hand side of the inequality (\ref{left-right})
with respect to \(\xs \in \partial \fs(0)\) obtain
\[
\partial \fs(0; \zk) \leq \zk p^k
\]
Assuming that \(\zk \to z^\bullet, p^k \to p^\bullet \) when \(k \to \infty\) and \(\gpk \to 0\)
according to Theorem \ref{thm-algo} obtain
\[
\partial \fs(0; z^\bullet) \leq z^\bullet p^\bullet.
\]
As \(p^\bullet \in \partial \fs(0) \)
by upper semi-continuity of the sub-differential mapping \(\partial \fs (\cdot)\)
\[
\partial \fs(0; z^\bullet) \leq z^\bullet p^\bullet \leq \sup_{p \in \partial \fs(0)} z^\bullet p = \partial \fs(0; z^\bullet)
\]
which implies that \( \zk p^k \to \partial \fs(0; z^\bullet) \) when \(k\to\infty\) or
\bleq{esti-theta}{
\partial \fs (\gpk; \zk) - \partial \fs(0; \zk) = \theta(\gpk, \zk) \to 0
}
when \(k \to \infty\).
Putting everything together we obtain
\[
\fs(0) - \xi_{k+1} \leq \theta(\gpk, \zk) \| g^k \|,\ \theta(\gpk, \zk) \to +0, \mbox{ when } k \to\infty.
\]
As \(x^k = -\gpk/(\fs(\gpk) - \xi_k) \) than due to upper semi-continuity of \(\partial \fs \) 
\bleq{gest}{
\|\gpk\| = \| x^k \| (\fs(\gpk) - \xi_k) \leq 2 \|X_\star\| (\fs(\gpk) - \xi_k) \leq  2 \|X_\star\| (\fs(0) - \xi_k)
}
and consequently
\[
\fs(0) - \xi_{k+1} \leq 2\theta(\gpk, \zk) \| X_\star\|(\fs(0) - \xi_k) = \lambda_k (\fs(0) - \xi_k)
\] 
with \(\lambda_k \to 0\) when \(k\to\infty\).
\EndProof

Next we consider the problem (\ref{copt}) with sup-quadratic objective function \(f\).
\begin{theorem}
Let 
objective function \(f\) in problem (\ref{copt}) is locally sup-quadratic with
sup-quadratic characteristic \(\tau\)
and \(\xi_k, k=1,2, \dots\) are defined by the Algorithm \ref{cco-algo}
with \(\xi_0 < -f_\star \).
Then \(\lim_{k\to\infty} \xi_k  = \fs(0) \) (Algorithm \ref{cco-algo} converges)
and for \(k\) large enough
\(\fs(0) - \xi_{k+1} \leq \tau^{-1}(\fs(0) - \xi_k)^2\)
(that is convergence is quadratic).
\label{thm-quad}
\end{theorem}
\Proof
It follows from local sup-quadratic behavior of \(f\) that \(\fs\) is differentiable in some neighborhood \(U\) of \(0\).
Therefore the subdifferentials of \(\fs(g)\) are singletons and we can consider \(\partial \fs(g) \) as just a vector.
It follows from sup-quadratic behavior of \(f\) 
that \( \| \partial\fs(g) - \partial\fs(0) \| \leq \tau^{-1}\|g\| \).
Consequently
\bleq{qud}{
\partial \fs (\gpk; \zk) - \partial \fs(0; \zk) = 
\partial\fs(\gpk) \zk  - \partial\fs(0) \zk \leq \| \partial\fs(\gpk) - \partial\fs(0) \| \leq \tau^{-1}\|\gpk\|
}
gives exactly \(\fs(0) - \xi_{k+1} \leq C (\fs(0) - \xi_k)^2 \)
with \(C = \tau^{-1}\).
\EndProof

Finally we consider the case of a sharp minimum in (\ref{copt}), namely that \( 0 \in \intr(\partial f(\xs)) \).
\begin{theorem}
Let the objective function of (\ref{copt}) has a sharp minimum at solution point \(\xs\),
all assumptions of the theorem \ref{thm-algo} are satisfied 
and \(\xi_k, k=1,2, \dots\) are defined by the Algorithm \ref{cco-algo}
with \(\xi_0 < -f_\star \).
Then there exists \(k^\star\) such that
\(\xi_{k^\star} = \fs(0) = -f_\star\).
\label{thm-fin}
\end{theorem}
\Proof
According to Lemma
\ref{lincon}
under conditions of sharp minimum
there is a neighborhood \(U\) of \( g = 0 \)
such that
\[
\fs(g) = \sup_x \{ gx - f(x)\} = g \xs - f(\xs)
\]
is a linear function of \(g\) in \(U\).

By Theorem \ref{thm-algo} \( \gpk \in U \) for \( k \) large enough
and let \(k^\star - 1\) is the first such index that \(g^{k^\star-1}_p \in U\).
Then
\[
\partial \fs(g^{k^\star-1}_p;g^{k^\star-1}_p) = g^{k^\star-1}_p \xs = \partial \fs(0;g^{k^\star-1}_p)
\]
and hence
\[
0 \leq \xi_{k^\star} - \fs(0) \leq
\partial \fs(g^{k^\star-1}_p;g^{k^\star-1}_p) - \partial \fs(0;g^{k^\star-1}_p) = 0
\]
and Algorithm \ref{cco-algo} terminates.
\EndProof
 
\makeatletter{}\section*{Conclusion}
\label{itogo}
The conceptual version of the dual epi-projection algorithm
has promising computational properties which makes it a viable candidate
for developing implementable versions.
First of all it guarantees global super-linear convergence to the optimum
for any solvable COP.
Second, it provides quadratic convergence and even finite termination
without any changes in the algorithm
for quite common types of COPs:
sup-quadratic, which strictly contain strongly convex,
and COPs with sharp minimum.
It is worth to notice that the algorithm is absolutely parameter-free,
use the first-order subgradient oracle only,
and requires no specific knowledge of any specific characteristics of COP,
like Lipshitz constants, strong convexity parameter or close enough
starting point.

The implementation perspectives for the algorithm depend upon the possibility
to produce practical version of the projection operator on \(\epifs\).
From the theoretical point of view it is easy to derive accuracy estimates for its termination
so it can be finitely solved for any required accuracy.
It can be used to preserve the overall rates of convergence
in terms of Algorithm \ref{cco-algo} iterations,
however the resulting computational complexity requires further investigations.

\end{document}